\documentclass[12pt]{article}
\usepackage{amsthm, amsmath, amssymb}
\newtheorem{theorem}{Theorem}[subsection]
\newtheorem{lemma}[theorem]{Lemma}
\newtheorem{corollary}[theorem]{Corollary}
\newtheorem{proposition}[theorem]{Proposition}
\newtheorem{remark}[theorem]{Remark}
\newtheorem{definition}[theorem]{Definition}

\newtheorem{hyp}[theorem]{Hypothesis}

\newcommand{\nc}{\newcommand}
\nc{\cH}{{\mathcal H}}
\nc{\cA}{{\mathcal A}}
\nc{\cG}{{\mathcal G}}
\nc{\cC}{{\mathcal C}}
\nc{\cO}{{\mathcal O}}
\nc{\cI}{{\mathcal I}}
\nc{\cB}{{\mathcal B}}
\nc{\cY}{{\mathcal Y}}
\nc{\cK}{{\mathcal K}}
\nc{\cX}{{\mathcal X}}
\nc{\cS}{{\mathcal S}}
\nc{\cE}{{\mathcal E}}
\nc{\cF}{{\mathcal F}}
\nc{\cZ}{{\mathcal Z}}
\nc{\cQ}{{\mathcal Q}}
\nc{\cN}{{\mathcal N}}
\nc{\cP}{{\mathcal P}}
\nc{\cL}{{\mathcal L}}
\nc{\cM}{{\mathcal M}}
\nc{\cT}{{\mathcal T}}
\nc{\cW}{{\mathcal W}}
\nc{\cU}{{\mathcal U}}
\nc{\cD}{{\mathcal D}}
\nc{\cJ}{{\mathcal J}}
\nc{\cV}{{\mathcal V}}
\nc{\bH}{{\mathbb H}}
\nc{\bA}{{\mathbb A}}
\nc{\bG}{{\mathbb G}}
\nc{\bC}{{\mathbb C}}
\nc{\bO}{{\mathbb O}}
\nc{\bI}{{\mathbb I}}
\nc{\bB}{{\mathbb B}}
\nc{\bY}{{\mathbb Y}}
\nc{\bK}{{\mathbb K}}
\nc{\bX}{{\mathbb X}}
\nc{\bS}{{\mathbb S}}
\nc{\bE}{{\mathbb E}}
\nc{\bF}{{\mathbb F}}
\nc{\bZ}{{\mathbb Z}}
\nc{\bQ}{{\mathbb Q}}
\nc{\bN}{{\mathbb N}}
\nc{\bP}{{\mathbb P}}
\nc{\bL}{{\mathbb L}}
\nc{\bM}{{\mathbb M}}
\nc{\bT}{{\mathbb T}}
\nc{\bW}{{\mathbb W}}
\nc{\bU}{{\mathbb U}}
\nc{\bD}{{\mathbb D}}
\nc{\bJ}{{\mathbb J}}
\nc{\bV}{{\mathbb V}}
\nc{\bbZ}{{\mathbb Z}}
\nc{\bR}{{\mathbb R}}
\nc{\fr}{{\rightarrow}}
\nc{\co}{{\overline{\nabla}}}

\begin{document}
                                %
                                %
                                %
\title{Monodromy of constant mean curvature \\surface
 in hyperbolic space
} %
\author{Gian Pietro Pirola
\footnote{Partially supported by
 1) PRIN 2005 {\em "Spazi di moduli e
 teorie di Lie"}; 2) INDAM;
3) Far 2006 (PV):
 {\em Variet\`{a} algebriche, calcolo
 algebrico,
 grafi orientati e topologici}.}
}
\date{}     
                                %
                                %
\maketitle
                                %
                                %
 \begin{abstract}
{\em
In this paper we give a global
version of the Bryant representation
of surfaces of constant mean curvature one
(cmc-$1$) in hyperbolic space.  This
allows to set the associated non-abelian
period problem in the framework of flat unitary
vector bundles on Riemann surfaces.
We use this machinery to prove
the existence of certain cmc-$1$
 surfaces having
prescribed global monodromy.} \vskip3mm
\vskip 1mm
 \noindent {\scriptsize {\bf Key words:}
.}
\end{abstract}
                                %
                                %

\thispagestyle{empty}           %
\pagestyle{myheadings}
\markright{Surfaces in hyperbolic space} %
                                %
                                %

 \subsection*{Introduction}

The local theory of
(cmc-$1$) in  hyperbolic space is
equivalent to the local theory of
minimal surfaces in  Euclidean space.
If this was originally the main
reason for their study
(but see also \cite{donaldson}),
the interest in the topic was
renewed by the fundamental work
of Robert Bryant \cite{bryant}.
The
  "Bryant-Weierstrass" formula
(see also \cite{lima} for the analysis of an earlier formulation)
allows to represent these
surfaces by holomorphic mappings.
To explain this, let
$SL(2,\mathbb C)$ be the special linear group
and $SU(2)$ the special unitary group.
We identify the quotient
$SL(2,\mathbb C)/SU(2)$
with the hyperbolic
$3-$space $\mathcal H^{3}.$
Letting $\pi:SL(2,\mathbb C)\to \mathcal H^{3}$
be the quotient
then \cite{bryant} a
simply connected cmc$-1$ in $ \mathcal H^{3}$
arises $S=f(U),$ where
 $U$ is an open set of the complex plane,
 $f=\pi\cdot g $ and
 $g:U\to SL(2,\mathbb C)$ is a holomorphic map.
In other words,
the entries of the matrix
$$
  g(z)= \begin{pmatrix} \nonumber
  a(z) & b(z) \\
  c(z) & d (z)
   \end{pmatrix},$$
 are holomorphic; moreover they satisfy
 the Bryant conditions:
\begin{equation} \label{one}
  \det g(z)=1,\ \ \det g'(z) =
  \det \begin
{pmatrix} a'(z) & b'(z) \\
  c'(z)  & d'(z)
   \end{pmatrix}=0 . \end{equation}

  The global theory of
 cmc-$1$ surfaces substantially differs from the
  theory of minimal surfaces in
 flat spaces. In fact,
 since $SU(2)$ is non abelian, the period problem
 cannot be solved by
 classical potential theory.
 In this paper
we  set the Bryant-Weierstrass representation
 formula in the framework of
unitary bundles on Riemann surfaces.
The period problem becomes then
equivalent to the existence
of suitable holomorphic sections
of a flat bundle. The holonomy
of the flat connection encodes the
monodromy data. To explain this,
let $X$ be a Riemann surface and
$f:X\to \mathcal H^{3}$ be
a conformal immersion such that
$f(X)$ is a cmc-$1$
surface. Then there are a rank $2$ vector bundle
$F$, with flat  $SU(2)$ connection
$\nabla,$ and
two holomorphic sections of $F$
$s_{1}$ and $s_2$ satisfying:
\begin{equation}\label{two}
\det(s_{1}, s_{2})= 1 \ ,
\ \det(\nabla s_{1},\nabla s_{2})= 0 .
\end{equation}
With respect to a unitary basis for the space
of harmonic sections of $F$,
the equations in (\ref {two}) become the
Bryant conditions (\ref {one}).
Fix a point $x\in X$ and let $\Pi_{1}(X,x)$ denote the
fundamental group of $X$ at $x$.
The monodromy map associated to
 $f$ is a homomorphism:
$$m(f):\Pi_{1}(X,x)\to SU(2),$$
which is the holonomy
of $\nabla.$
Conversely, two holomorphic
sections of a flat $SU(2)$-bundle on $X$
that satisfy the equations (\ref{two})
define a cmc-$1$ surface
in
 $\mathcal H^3.$
\medskip

In the case of an algebraic
complex curve, the flat unitary bundles are
equivalent to certain stable parabolic holomorphic vector
bundles. This allows, via
a Riemann-Roch type theorem,
to prove the existence
of surfaces with fixed monodromy. More precisely,
let $Y$ be a compact Riemann surface and
 $H\subset Y$ be any finite
subset, possibly empty.
The  result  given in \ref{esistenza} implies
(see the end of section 4)
the following:
\bigskip

\noindent {\bf Existence Theorem}\
 \noindent {\em Let $m:\Pi_{1}(Y\setminus H) \to SU(2)$
 be a group homomorphism.
  Assume that the
image of $m$ is not an abelian group. Then there
 is a finite set $D\subset Y$
and an immersion
$f: Y\setminus (D\cup H)\to  \mathcal H^{3}$
such that
 $f(Y\setminus (D\cup H))$ is a
cmc-$1$ surface and
$m(f):\Pi_{1}(Y\setminus (D\cup H))
\to SU(2)$ factors through $m.$
That is $m(f)= m j_{\ast}$ where
$$j_{\ast}:
\Pi_{1}(Y\setminus (D\cup H))\to
\Pi_{1}(Y\setminus  H)$$
 is induced by the inclusion.}
 \medskip

Moreover in \ref{esistenza}
 the sections  are meromorphic
at
 $D\cup H,$ the {\em ends}
 of $f.$
The study of
meromorphic sections satisfying
the Bryant conditions (\ref{two})
turns out to be
quite interesting.
Our proof of the existence result
relies on a homogenization
procedure.
This allows to use intersection theory
in projective spaces.
We note we are not able
to prove that our surfaces are
complete: if we complete them,
we are no longer able to prove
they are immersed.
This problem should be clarified by an analysis
of the moduli space parametrizing
sections that satisfy
Bryant conditions. In contrast with the
Euclidean minimal surfaces  case (see \cite{pi}), an
infinitesimal useful description has not yet been found.
\medskip

The paper is organized as follows.
In section $1$, starting from simple consideration
on complex Lie groups,
we obtain the global version
of the Bryant-Weierstrass
representation \ref{global}.
In section $2$ we analyze the case of a punctured disk
and the induced local parabolic structure.
In  \ref{polo} we give a simple, but important, remark
on the Laurent series of our sections.
In section  $3$ we recall the
correspondence between stable and
flat $ SU(2)$ bundle on a smooth complex
algebraic curve.
We remark the result of lemma \ref {propria}. It provides
a properness result in the case of irreducible connections.
In section $4$, all previous results are collected to prove
our main Theorem \ref{esistenza}.
A more general setting is considered in
section $5,$ which is in certain sense
a complementary section.
We use the unifying language
of Higgs fields to consider
both the surfaces introduced in \cite{ktuy}
and  the periodic cmc-$1$ case.
This allows to formulate
a general non-abelian period problem
(see \ref{higgs}).

\medskip

It is a great pleasure to thank Maasaki
Umehara for his precious
advice. I am really grateful to him.
I would like to thank Enrico Schlesinger, 
who read the manuscript and suggested
many improvements.

\section {Bryant representation}
                                %
                                %
                                %
\subsection{Lie groups and flat structures}
Let $X$ be a smooth connected Riemann surface
and   $\mathcal O_{X}$
be the holomorphic structure sheaf of $X$.
Let $G$ be a complex Lie group and let
${\it g}$ be its Lie algebra.
We denote by $J:{\it g}\to {\it g}$
the map induced by the complex structure.
Let $H$ be a closed (Lie)
subgroup of $G$ with Lie algebra ${\it h}$.
 We shall assume $H$ transverse to the complex structure
 of $G$, that is:
 \begin{equation} \label{transverse}
  J{\it h}\cap {\it h} = 0_{{\it g}}.
   \end{equation}
Letting $\pi : G\to G/H$ be the quotient,
 we consider maps $f: X\to G/H.$ We say that $f$
 is $h$-liftable
  if there are local holomorphic
  liftings of $f$ to $G.$ This means that
  there exist an open covering
  $\{ U_{\alpha} \} _{\alpha\in I} $
  of
  $X$ and holomorphic maps:
\begin{equation}
\label{lift}
 g_{\alpha}: U_{\alpha} \rightarrow G
\end{equation}
  such that $\pi g_{\alpha}= f|U_{\alpha}.$
  For $\alpha$ and $\beta \in I$ we define
 \begin{equation} \label{cociclo}
  g_{\alpha \beta}=g_{\alpha}^{-1} \cdot
 g_{\beta} : U_{\alpha}\cap
  U_{\beta} \rightarrow H\subset G.
 \end{equation}
  Since $g_{\alpha \beta}$ is holomorphic
   and $H$ is
 transverse to the complex structure,
 the $g_{\alpha \beta}$ are
 locally constant and therefore define a
 $H$-principal flat bundle on $X.$
 If  a point $p$ of $X$  is fixed, loops restriction
 induces
 the monodromy homomorphism :
 $$
 m : \Pi(X,p) \rightarrow H\subset G.
 $$
If $V$ is a complex vector bundle and
  $\rho: G\to GL(V,\mathbb C)$ is a complex
  linear representation, the
  $\rho(g_{\alpha \beta})$ are transition
  functions for a flat
  complex vector bundle $F$ on $X.$
  The structure group of $F$ reduces
  to $\rho(H).$ Then $F$ has a natural holomorphic structure
  and the
  composition
 \begin{equation} \label{rep}
  \theta=\rho\cdot m : \Pi(X,p) \to \rho(H)\subset
 GL(V)
 \end {equation}
  is the associated monodromy representation.


\subsection {Special and unitary group}

Let $M(n)$ denote the vector space of square
complex matrices of order $n$.
In the sequel, unless specified differently,
 $G$ and  $H$ will denote
 respectively  the special
linear group $G=SL(n,\mathbb C)=
\{g\in M(n) : \det(g)=1\}$,
and the special unitary group
 $H=SU(n)=\{g \in SL(n,\mathbb C): g\cdot
g^{\ast}=e\}$, where
$e\in M(n)$ is the identity
matrix and  $g^{\ast}=\!\!^{t}\bar{g}$
is the adjoint of $g$. Clearly,
$H$ is transverse to the complex structure of $G$: $h=su(n)$.
By the spectral theorem we identify the
quotient $SL(n,\mathbb C)/SU(n)$ with
$$L(n)=\{g\in SL(n,\mathbb C): g=g^{\ast}, g>0\},$$
so that $$ \pi(g)= g\cdot g^{\ast}.$$

 Let $f:X \rightarrow L(n)$ be a
  $h$-liftable map and let
 $\{g_{\alpha \beta}\}$ be the
 $ SU(n)$-cocycle defined
 in (\ref {cociclo}).
Let $W= \mathbb
C^{n}$ be the standard
representation of $SU(n)$ and $F$ be
the associated
 flat complex vector bundle on $X.$ Note
that $F$ is a rank $n$ holomorphic
vector bundle on $X$ with trivial
determinant. The natural hermitian
product on $F$ will be denoted by $<,>.$
For a fixed $v\in W $,
the collection
$\{g_{\beta}^{-1}v\}_{\beta \in I}$
 (see formula \ref{lift})
defines a holomorphic section of
$F$. Letting $H^0(X,F)$ denote the vector
space of the holomorphic
sections of $F,$ we thus obtain an
inclusion $j:W\to H^0(X,F).$ We let
$\{v_i\}_{i=1,\ldots, n}$ be the canonical
basis of $W$, and set
$j(v_{i})= e_i.$ Under the isomorphism
$\det(F)= \mathcal O_{X}$ we
have:
\begin{equation} \label{wedge}
 e_{1}\wedge \ldots \wedge e_{n}=1.
\end{equation}
If we denote by $F_{x}$ the fiber of $F$
over $x\in X$ the $e_{i}(x)
$ give a $SL(n,\mathbb C$) basis of $F_{x}$
at every point $x\in X.$
\begin{definition}
\label{specialframe}
 A holomorphic frame $\{e_{i}\}_{i=1,\ldots, n}$
 satisfying (\ref{wedge}) will be called a
 special frame of $F$.
  \end{definition}

 To any $h$-liftable map
 $f:X\to L(n)= Sl(n,\mathbb C)/ SU(n)$
 is associated
 a special frame $ e_{1} \ldots e_{n}$
 of a $SU(n)$  flat vector bundle
 $F$. Conversely, if $F$ is a $SU(n,\mathbb C)$
 flat vector bundle on
 $X$ and $ e_{1} \ldots e_{n}$
is a special frame of $F,$ for every $x\in X$ we define
 $$
 A(x)= (e_{1}(x),\ldots, e_{n}(x))\in
 Hom(\mathbb C^{n},F_{x}).
 $$
This gives a holomorphic section
  $A\in Hom(\mathbb C^{n},F)$. On the other hand,
the adjoint $ A^{\ast}(x):$ 
$$<A(x)v,w>_{x}=<v,A^{\ast}(x)w>,$$
provides an
antiholomorphic section $A^{\ast}$ of 
$Hom(F,\mathbb C^{n}).$ Then
$$f(x)=
 A(x) \cdot A^{\ast}(x)
 \in L(n)=Sl(n,\mathbb C)/ SU(n)$$
 defines a
$h$-liftable map
 $f:X\to Sl(n,\mathbb C)/ SU(n),$
the holomorphic
local lifting of $f$ being  $A(x)$
 written in local coordinates.

\begin{definition} \label{equi}
 Two special frames
$\{e_{i}\}_{i=1,\ldots, n}$ and
$\{f_{i}\}_{i=1,\ldots, n}$ of $F$ will
be said equivalent if there is
$U\in SU(n)$ such that $U{e_{i}=f_i \ \
i=1\ldots, n}.$
 \end{definition}
Special equivalent frames
correspond to the same $h$-liftable mapping
 $f,$ thus we have:

\begin{proposition}
 \label {lol}
 There is a one to one correspondence
 between $h$-liftable maps
 $$f:X\to SL(n,\mathbb C)/ SU(n)$$
 and equivalence classes of special
 frames of $SU(n)$-flat vector bundles of $X$.
  \end{proposition}
  \begin{proof} Clear.
  \end{proof}

 The monodromy (\ref {rep}) gives
  a map
 $\theta : \Pi(X,p) \to SU(n).$ The group
 \begin {equation} \label{holonomy}
 G_{\theta}= im(\theta)
 \end{equation}
is well defined up to conjugation
 and determines the flat bundle $F.$

 \begin{definition}
 \label{rid}
 The flat bundle $F$ is said reducible
 (resp.  irreducible) if
 the corresponding representation  $\theta $
 of the fundamental group is
 reducible (resp.  irreducible).
  \end{definition}

 \begin{remark} \label{periodic}
There are many Lie subgroups
$H \subset SL(n,\mathbb C)$
transverse to
the complex structure,
for instance $SL(n,\mathbb R)$. For our
purposes, however, it is useful to assume
$h= su(n)$, that is, the connected
component $H_e$ of $H$ is the special
unitary group: $H_{e}=SU(n)$. In
this case $K=H/H_{e}$ is a discrete
group acting on $L(n).$ We have
$L(n)/K=G/H.$ If $f:X\to G/H$
is as before and $K$ acts freely on
$L(n)$, by taking a suitable covering
$\tilde{X}\to X,$ we may define
 a  periodic lifting (see also section 5)
$\tilde{f}: \tilde{X} \to
L(n)$.
\end{remark}

\subsection {Flat connections}

We fix a flat $SU(n)$ vector bundle $F$ on
$X$ with hermitian metric $<,>.$ Denote by
$C^{\infty}(F)$ the sheaf of smooth sections of $F.$ Let
$A^{k}(F)=C^{\infty}(F\otimes\Omega^{k}_{X})$
be the sheaf of $k$-forms on $F$.
Using the complex structure on $X$ we can also consider the space
$A^{p,q}(F)$ of $p,q$ forms on $X$ (\cite {griffithsha}).
  Let $\nabla: A^{k}(F)\to A^{k+1}(F) $
(see \cite {Kobayashi}) be the associated
unitary flat connection  $$\nabla^2=0,$$
 which preserves the metric
$$ d <s,t>=<\nabla s,t>+ <s, \nabla t>$$
and the complex structure:
$$
\nabla''= \bar{\partial}.
$$
Here $\nabla''$ denotes
the composition of $\nabla:
C^{\infty}(F)\to A^{1}(F)$ and
the projection $A^{1}(F)\to A^{0,1}(F).$
Let $$\ker \nabla \subset C^{\infty}(F)
$$
be the sheaf of harmonic sections of
 $F$. We have $\ker \nabla\subset
\ker \nabla''.$
We  identify as usual
$ F $ with the sheaf  $\ker \nabla''$
of  holomorphic sections of $F$.
\begin {remark} \label{ridu2}
The holonomy of $\nabla$  associated to $F$
  is (up conjugation)
the group $G_{\theta}$ (see \ref{holonomy}).
 We say that $\nabla$ is
 irreducible if the monodromy
 representation  \ref{rid}
is irreducible.
\end{remark}



\subsection {Bryant condition}

We now set $G=SL(2,\mathbb C)$
and $H=SU(2).$
 The quotient $L(2)=G/H$ can
be identified with the hyperbolic
$3-$space $\mathcal H^{3}.$ By \ref{lol}
$h$-liftable maps $f: X\to \mathcal H^{3}$
correspond to special
frames of $SU(2)$ vector bundles on $X.$
We fix a flat $SU(2)$ vector bundle $F$ on
$X$, and denote by  $\nabla$ its connection.
We remark that
 $\nabla$, (see \ref {ridu2}),
is reducible if and only if $G_{\theta}$
is abelian and
hence, up to conjugation, $G_{\theta}$ is
contained in
the subgroup of diagonal matrices.
\medskip

If $s$
is a holomorphic section of $F,$
$s\in H^{0}(X,F),$
 then $\nabla s$ is a holomorphic
 section of $F\otimes \omega_{X},$  $\nabla
 s\in H^{0}(X,F\otimes \omega_{X}).$
 Here $\omega_{X}$ denotes the sheaf
 of holomorphic differentials on $X.$
 Since $\det F $ is trivial,
 $\det(F\otimes \omega_{X})=\omega_{X}^{2}$
 is the sheaf of holomorphic
 quadratic differentials of $X.$
If $s\in H^{0}(X,F)$ and $t \in
 H^{0}(X,F)$ we define:
\begin{equation}
\label{quadratic}
 \Omega (s,t)= \nabla  s\wedge \nabla t \in
H^{0}(X,\omega_{X}^{2}).
\end{equation}
The vanishing of the above differential
(\ref{quadratic}) has an important
 geometric meaning (see \cite{bryant}).
 Recall that a special frame (see
 \ref{specialframe}) of $F$ is a couple
 of holomorphic sections $s$ and
 $t$ such that
$$s\wedge t=1.$$
\begin{definition}
 \label {bryant}
  We say the special frame $(s,t)$
   of $F$ satisfies the Bryant
  condition if
  $ \nabla s\wedge \nabla t = 0.  $
In this case $(s,t)$ will be called a
$\mathcal B$-frame.
\end{definition}
Let $(s,t)$ be $\mathcal B$-frame of $F$.
Fix a point $\bar x\in X$ and let $\{U(\bar x),z\}$
be a simply connected
open coordinate neighborhood of $\bar x$. Let
$z:U(\bar x)\to \mathbb C$ be the coordinate map.
We can find, by
parallel transport, a harmonic
unitary frame $h_{1}$, $h_2$ of $F$,
$\nabla h_i=0$, on
$U(\bar x)$ such that $h_{1}\wedge h_{2}=1.$
We may write:
$$ s=a(z)h_1 + c(z)h_1; \quad t = b(z)h_1 + d(z)h_1 .
$$
hence
 $$ \nabla s =(a'(z)h_1 + c'(z)h_2)dz ; \quad \nabla
t = (b'(z)h_1 + d'(z)h_2)dz .$$
Here $a,b,c$ and $d$ are
holomorphic functions. Define
the matrix
$$  A(z)= \begin{pmatrix}
  a(z) & b(z) \\
  c(z) & d (z)
   \end{pmatrix},$$
and note $\det A(z)=1.$ Taking derivative
we may write the Bryant condition as
  \begin {equation}
  \det A'(z) = \det \begin
{pmatrix} a'(z) & b'(z) \\
  c'(z)  & d'(z)
   \end{pmatrix}=0 \label{brlocale}. \end{equation}
 Let
 $f: U(x) \to \mathcal H^3$
(see \ref{lol})
 $$f(z)=
A(z)\cdot A^{\ast}(z),$$ be the associated
map.
 Assume that $f$ is non-constant, then
following \cite {bryant} (see also
 \cite{small}, \cite{uy} and the H.
  Rosenberg
contribution in \cite{Meeksrosrosenberg}),
we observe that
(\ref {brlocale}) gives the Bryant-Weierstrass
representation of the
constant mean curvature one (cmc-$1$) surface
$f(X)\subset \mathcal H^{3}.$ We have proven:
  \begin{proposition}
   The maps  $ f: X\to \mathcal H^{3}$ such that
  $f(X) $ is a (branched) cmc-$1$ are
  in one to one correspondence with
  equivalence classes (see \ref{equi}) of non-trivial
   $\mathcal B$-frames
  \ref{bryant}.
  \end{proposition}
  One verifies that, if $(s,t)$ is a $\mathcal B$-frame,
  then $f:X\to
  \mathcal H^{3}$ is an immersion if and only
  if all the local holomorphic
  lifting
  $ g_{\alpha}: U_{\alpha} \rightarrow SL(2,\mathbb C)$
  are
  immersions (see \cite{Meeksrosrosenberg}).
  It follows
  then that $x$ is a branch point if and only
  if $A'(z(x))=0$ (see \ref{brlocale}).
   Therefore  $f$ is an immersion if and only if
    the set of branch points of $f$
  $$Z(f)= \{x\in X:
 \nabla s(x)=\nabla t(x)=0\}$$
 is empty. We
 finally state the following
 (compare with \cite{bryant}):

  \begin{proposition} \label{global}
  {\bf Global Bryant-Weierstrass
  representation} Let $M \subset \mathcal H^3$
  be a cmc-$1$ immersed
  surface. Then $M=f(X)$ where $X$
  is a Riemann surface, and $f:X\to
  \mathcal H^3$ is the $h$-liftable
  map  associated to a couple of sections
  $s,t$ of a $SU(2)$ flat bundle $F$ such that
   \begin {enumerate}
  \item [a)]
  $ {\bar \partial s}= {\bar \partial t}=0$
  (holomorphicity); \item[b)] $s\wedge t=1$
  (special-frame);
  \item[c)]
  $\nabla s\wedge \nabla t=0$
  (Bryant condition);
  \item[d)]
  $Z=\{x\in X: \nabla s(x)=
  \nabla t(x)=0\}= \emptyset$ (immersion).
 \end{enumerate}
  \end{proposition}



\section{The punctured disk: the end of a surface}\label{disk}

Let $\Delta =\{ z\in \mathbb C:
|z < 1\}$ be the unit disk and  denote by $O$ the origin of $\Delta$. Set
$\Delta^{\ast}= \Delta - \{O\}.$
We will recall the classical description of
 the $SU(2)$ connections of
 $\Delta^{\ast}$ which will be used in our setting.

\subsection{Singular connection on  $\Delta$}

Since the fundamental group of $\Delta^{\ast}$
is cyclic, any $SU(2)$ bundle over $\Delta^{\ast}$ is reducible (see \ref{rid}
 and \ref {ridu2}) and the holonomy map can
be described by a single matrix
$T(a):$ $$ T(\alpha)=
\begin{pmatrix}
\exp(\pi i\alpha) & 0\\
0 & \exp(-\pi i\alpha)
\end{pmatrix}
$$ where $\alpha \in \mathbb R,\ 0\leq \alpha <1.$
The holonomy group is
$G_{\alpha}= \{ T(\alpha)^n: n\in \mathbb Z \}.$
To construct a flat bundle with
$G_{\alpha}-$holonomy, take the trivial bundle
 $$\mathbb C^2 \times \Delta \to \Delta.$$
Then we define the singular connection:
\begin{equation} \label{consing}
{\overline\nabla}_{\alpha} = d-
\begin{pmatrix}
 \alpha & 0 \\
 0 & -\alpha
  \end{pmatrix}\frac{dz}{z},
 \end{equation}
 that is,
 $$
 {\overline\nabla}_\alpha(f,g)=(df-\alpha \frac{f}{z}dz,
dg+\alpha
\frac{g}{z}dz).
$$
Let  $\mathcal O_{\Delta}(nO)$  be the
  sheaf of meromorphic function on $\Delta$ that are
   holomorphic on
  $\Delta^{\ast}$
  and  have at most at most a pole
  of order $n$ at the origin $O$ of $\Delta$.
  We have $$\overline\nabla: \mathcal O_{\Delta}^{2} \to
 \mathcal O^{2}_{\Delta}(O)\otimes \omega_{\Delta}.
 $$
Its residue
 \begin{equation} \label{res}
res{\overline\nabla}_{\alpha} = \Gamma_{\alpha}=
\begin{pmatrix}
 \alpha & 0 \\
 0 & -\alpha
 \end{pmatrix}
 \end{equation}
 is an operator on  $\mathbb C^2,$
  the fiber of $O,$ having eigenvalues
$\alpha$ and $-\alpha$, and, if $0<\alpha<1 $,
eigenspaces $L_{\alpha}=\{(z,0): z\in \mathbb
C \}$ and  $L_{-\alpha}=\{(0,z): z\in \mathbb C \}.$
\medskip
 Let $V=\mathbb C^2|_{\Delta^{\ast}}$
be the restriction on $\Delta^{\ast}$  we
set $$\nabla_{\alpha} =
{\overline\nabla}_\alpha|_{\Delta^{\ast}}.$$
Note that
 $\nabla_{\alpha}$ is a flat connection on $V$.
The matrix
$$
H(\alpha) =
\begin{pmatrix}
w^{\alpha} & 0 \\
0 & w^{-\alpha}
\end{pmatrix}
\ w= z\cdot {\bar z}= x^{2}+y^{2},\ z=x+iy,
$$
 defines an hermitian product $<,>_{\alpha}$ on $V$
 compatible with
  $\nabla_{\alpha}.$ The holonomy group of
  $\nabla_{\alpha}$ is
  $G_{\alpha}.$ We remark that the metric
  $<,>_{\alpha}$ extends
  continuously by $0$ on $L_{\alpha}$
  and gives a natural filtration
  $$ 0\subsetneq L_{\alpha} \subsetneq \mathbb C^2_{O}.$$
We have:
\begin {proposition}
The connection $\nabla_{\alpha}$ is a $SU(2)$
flat connection with holonomy
$G_{\alpha}.$
\end{proposition}
\begin{proof}
One has  $\nabla_{\alpha}^{2}=0,$
$\nabla''_{\alpha}={\bar\partial}$
and
 $\nabla_{\alpha}$ is compatible
 with  $<,>_{\alpha}.$
\end{proof}
\medskip

\begin {remark}
Set
$\widetilde{d} =\wedge^2{\overline\nabla}$. Then $\widetilde{d}$
defines the trivial
connection on $\det E.$ In fact, let
$v_{1}=(0,1)$, $v_{2}=(1,0)$ be the
 constant sections, $u=v_1\wedge v_2$ and fix a
function $b$. We have
$$\widetilde{d}(b\cdot u)=
d (b (v_1\wedge v_2)) + {\overline\nabla}_{\alpha}(v_{1}) \wedge v_{2}+
v_{1}\wedge {\overline\nabla}_{\alpha} (v_{2})=$$
$$(db)\cdot (v_1\wedge v_2) +
\frac{\alpha dz}{z}(-v_{1}\wedge v_{2}+v_{1}
\wedge v_{2})=(db)\cdot u,$$
that is, $\widetilde{d}=d.$
\label{bandet}
\end{remark}
If $\{F,\nabla\}$ is a $SU(2)$ flat
 bundle on $\Delta^{\ast}$ then it is
isomorphic to $\{ V, \nabla_{\alpha}\}$
for some $\alpha,$
$0\leq \alpha<1.$
Therefore there is a natural logarithmic extension
$\{E,{\overline\nabla}\}$
of $\{F,\nabla\}$ to $\Delta$:
 $$
 \overline{\nabla}:
E \to E(O)\otimes \omega_{\Delta},
$$
where $E$ is isomorphic to the trivial
bundle $\mathbb C^{2}$ and
$E(O)= E\otimes \mathcal O_{\Delta}(O).$
Moreover, if $E_{O}$ is the fiber
of $E$ at $O$,
$\Gamma=res \overline{\nabla}\in Hom(E_{0},E_{0})$
 has eigenvalue
$\alpha$ and $-\alpha.$
If the holonomy is non-trivial, i.e.
$\alpha>0$,   let
 $W=\ker( \Gamma-\alpha
Id)$ where $id$ is the identity.
The hermitian product on
$E$ extends by $0$ on all of $W$.
We have a filtration
$E_{0}\supsetneq W\supsetneq 0.$

\begin{definition} \label{locpar}
Suppose $\{F,\nabla\}$ is
 a $SU(2)$ flat bundle $\Delta^{\ast}$ with
 extension $\{E,{\overline\nabla}\}
 \simeq \{ \mathbb C^{2},{\overline
 \nabla_{\alpha}} \}$,  $0 \leq \alpha <1.$
 If $\{F,\nabla\}$  has non
 trivial holonomy, i.e. $\alpha>0$,
 the associated local
 parabolic structure on $E$ is given by:
\begin{enumerate}
\item[i)] the weight \ $\alpha\in ]0,1[\subset \mathbb R ;$
\item[ii)] the filtration
 $E_{0}\supsetneq W\supset 0$ \ where
 \ $W=\ker(\Gamma-\alpha\ id)$ \
 and \
 $\Gamma = \ res {\overline\nabla}_{\alpha}$.
\end{enumerate}
\end{definition}


 \subsection {Bryant condition on $\Delta^{\ast}$}

We write  Bryant condition for
sections of $\{V,\nabla_{\alpha}\}$. Let
$s$ and $t$ be holomorphic sections of $V$.
Laurent expansion allows
to write $s =(a(z), b(z))$ and $t= (c(z), d(z))$
as singular sections
of $\mathbb C^{2}.$ Assume $s \wedge t=1,$ i.e.
$ad-bc=1.$ Let
$f:\Delta^{\ast}\to \mathcal H^{3}$ be the
 $h$-liftable map associated
to this special frame \ref{lol}.

 \begin{definition}
 We say  the frame $(s,t)$
 has finite type if $s$ and $t$
extend meromorphically to $\Delta$.
  The order of $(s,t)$ is the smallest
integer $ n\geq 0$
such that $z^{n}s$ and $z^{n}t$ are holomorphic.
We also say
the $h$-liftable associated map (or cmc-$1$
immersion) has order
$n$.
\label{mero}
\end{definition}

\noindent Denote by $ \Omega_{\alpha}(s,t)
$ - compare formula  (\ref {quadratic}) above - the singular extension of
$\nabla_{\alpha} s
\wedge \nabla_{\alpha} t$ to $\Delta$:
 \begin{equation}
\Omega_{\alpha}(s,t)=
\overline\nabla_{\alpha} s\wedge \overline\nabla_{\alpha}t=
\det \Big[\begin
{pmatrix} a' & b' \\
c'  & d'
\end{pmatrix} -  \alpha \frac{1}{z} \begin {pmatrix}
a & -b   \\
c & -d \end{pmatrix} \Big](dz)^{2}.
\label{omegadisco}
 \end{equation}

Let
$$H^0(\Delta,\omega^{2}_{\Delta}(m))=
H^0(\Delta,\omega^{2}_{\Delta}\otimes \mathcal O_{\Delta}(mO))$$ be
the
space of holomorphic quadratic differentials
on $\Delta$ having a pole
of order at most $m$ at $O.$
If $(s,t)$ has finite type of order $n$,
 then $
\Omega_{\alpha}(s,t) = \overline
 \nabla_{\alpha} s\wedge \overline \nabla_{\alpha}t \in
H^0(\Delta,\omega^{2}_{\Delta}(2n+2)).$

 The following more precise result will play a role in
our existence result:

 \begin{lemma} \label{polo}
 Assume  $s, t\in H^0(\Delta,\mathcal
O_{\Delta}(nO))$,
and $s\wedge t=1$ (or else $s\wedge t$ holomorphic at $O$).
We have
 \begin{enumerate}
 \item   $ \Omega_{\alpha}
 (s,t)\in H^0(\Delta,\omega^{2}_{\Delta}(2n+1))$
  if $0<\alpha< 1 ;$
   \item $ \Omega_{0} (s,t)\in
  H^0(\Delta,\omega^{2}_{\Delta}(2n)).$
  \end{enumerate}
 \end{lemma}
 \begin{proof} Write the condition
that $s\wedge t$ is holomorphic
  on the Laurent series of $s$ and $t.$ Then compare
  with the Laurent series of
   equation  \ref{omegadisco}.
  \end{proof}

 \medskip Now assume that $(s,t)$ is a Bryant
frame, that is,  $s\wedge t=1$
and  $ \Omega_{\alpha} (s,t)=0$ \ref{bryant}.
We obtain the: \medskip

\noindent {\bf Bryant equations for the punctured
disk:} \begin {equation} \left \{
 \begin{array}{ll}
ad-bc=1 \\
z^{2}(a'd'-b'c') +z \alpha (ad'-a'd-b'c+bc') - \alpha^{2}=0.
\end{array}
\right .  \label{stareq}
\end{equation}


\section{Flat bundles on algebraic complex curves}

Non-abelian Hodge-theory
gives correspondences between
representations of the fundamental
group of algebraic (or compact K\"ahler)
varieties, that is, flat vector bundles,
and algebraic objects: (parabolic)
stable (Higgs) bundles.
We will only describe the $SU(2)$ case on algebraic curves.
\subsection{ Parabolic and stable bundles}

Let $Y$ be a compact Riemann surface.
A real divisor of $Y$ is a finite
 combination $S=\sum a_{i}P_{i},$ $a_i\in
\mathbb R,$\  $P_{i}\in Y.$
We will always assume that $P_{i}\neq P_{j}$ if
$i\neq j.$ We say that $S$
is effective if $a_{i}\geq 0$ for all $i.$
The
support of $S$ is  the divisor
$$supp \,(S) = \sum P_{i}: a_{i}\neq
0.$$ A divisor $D$ is simple if
 $D=supp \,(D).$ We  identify as usual
the set $ \cup_{i=1}^{n} \{P_i\}$ with the simple divisor
$D= \sum P_{i}$
(niding the case $D=0$). If $D$ is simple,
$X=Y-D$ is an algebraic complex curve.  \medskip

 Let $E$ be a
holomorphic vector bundle on $Y$.
 Fix a simple divisor $D=\sum
P_{i},$ let $E_P$ be the fiber of
$E$ at $P\in Y.$ Set $E_i=
E_{P_{i}}$ if $P_{i}\in D.$ We shall consider only the parabolic
structures, which we call {\em special}, connected with
 the $SU(2)-$bundles.
 For the general case
we suggest \cite {methases} and also \cite {belkale}.

 \begin {definition} \label{parabolic}
 A special parabolic vector bundle
 is given by the data
  $ \{ E, \mathcal P\} \equiv
  \{ E, D(\mathcal P), W_{i}\},$
  where:
 \begin{enumerate}
  \item[a)] $E$ is a rank
  $2$ holomorphic vector bundle such that $\det
  E=\cO_{Y};$
  \item[b)] $D(\mathcal
  P)=\sum_{i}\alpha_{i}P_{i}$ is a real weight
  divisor where $0 < \alpha_{i}< 1 ; $
  \item[c)]  $W_{i}$ is a
  proper subspace of $E_i$,  $0\subset W_i\subset E_i$, for every $i$.
   \end{enumerate}
  We set $D_{\cP}=supp \, D(\mathcal P)=\sum_{i} P_{i}$ and 
   $d_{\mathcal P}=\sum_i \alpha_i.$
  The trivial parabolic structure on $E$ is
given by the divisor $D = 0$.
\end{definition}

Since in our case the weight sequence at any point has just one element (see \cite {belkale}),
we have encoded this information into a real divisor.
We would like to stress that this is possible only in the rank $2$ case.
Let $ L\subset E $ be a
holomorphic sub-line bundle.
Let $L_{i}\subset E_{i}$ be the fiber of $L$ at $P_i\in D.$
We set
\begin{eqnarray}
      \left\{ \begin{array}{ll}
      \gamma_i=\ \alpha_i \quad \ if \quad L_i=W_i \\
      \gamma_i=-\alpha_i \quad   if \quad   L_i\neq W_i.
      \end{array}
    \right. \label{parL}
  \end{eqnarray}
and define the real divisor
 $$
 D(L,\mathcal P)= \frac {1}{2}\sum_{i=1}^{n}
\gamma_i P_{i},
$$
which has degree  $\deg D(L,\mathcal P)=
\frac {1}{2}\sum_{i=1}^{n} \gamma_i$.
We define the parabolic degree $par(L)$ of $L$ by the formula
 \begin{equation}
  \label{pardeg} par(L)= \deg L + \deg D(L,\mathcal P).
  \end{equation}
We give the following:

\begin{definition}
A special parabolic bundle $\{E,\mathcal P\}$
 is called stable (respectively
semistable)  if
 $par(L)<0$ (resp.  $par(L)\leq 0)$ for every holomorphic line bundle
 $L\subset E.$
\end{definition}

\begin{remark}

With the trivial parabolic
 structure on $E$, we have  $par(L)= \deg L$.
Since $\det E=\cO_{Y}$,  parabolic
stability in this case coincides with usual stability, see
\cite{methases}.
\end{remark}

   We prove now a standard
result we will use in section $4$.
  For any holomorphic
 vector bundle $\cK$ on $Y$ we set
$h^{i}(\cK)= \dim H^{i}(Y,\cK).$

\begin{proposition} \label {vanishing}
Let $\{E,\mathcal P\}$ be a special stable parabolic bundle and
$L$ be a line bundle of
degree $d$.
\begin{enumerate}
\item If $d\geq\frac{1}{2}d_{\mathcal P} +2g-1$,
then  $ h^{1}(E(L))=0.$
\item If $d\geq \frac{1}{2}d_{\mathcal P} +2g$,
then $E(L)$
is generated by its global sections.
\end{enumerate}
\end{proposition}
\begin{proof} \

\begin{enumerate}
\item By Serre duality we have to prove that
$ h^{0}(E(L^{-1})\otimes\omega_{Y})=0.$
Note that  $\det E=\cO_{Y}$
and hence that
$E$ is self-dual. Assume by
contradiction that $s\neq 0$ is a
global section of
$E(L^{-1})\otimes\omega_{Y}.$ It
 defines a map
 $ L\otimes \omega^{-1}_{Y} \to E.$
  The parabolic degree
of the image line bundle is bigger than
$\deg L- \deg\omega_{Y} - \frac{1}{2}
\deg D(L,\mathcal P)\geq 0$. This gives a contradiction.
\item Fix a point $P\in Y.$
 Consider the exact sequence
$0 \to
E(L(-P))\to E(L) \to E(L)_{P}\to 0.$
We have $H^{1}(E(L(-P)))=0$
 by the first part of the proof.
Hence the map $ H^{0}(Y,E(L))\to E(L)_{P}$
is surjective. \end{enumerate} \end{proof}



\subsection{Flat bundles and parabolic structures}

Let $X=Y-D$, where $Y$ is compact and
$D=\sum P_{i}$ is simple. Let
$\{F,\nabla\}$ be an $SU(2)$ flat vector bundle
on $X,$ $$\nabla: F \to F\otimes \omega_X$$
is flat:
 $\nabla^2=0$ and there is an hermitian
 form $<,>,$ on $F$ compatible with
$\nabla.$
We use the local
construction described
in section \ref{disk}.
Take a coordinate disk $ U_{j}$
around a point $P_{j},$ then
consider $U^{\ast}=U_{j}-\{P_{j}\}.$ The restriction of
$\{F, \nabla \}$
to  $U_{j}^{\ast}$
extends (see \ref
{locpar}) to a singular connection defined on $U_{j}$.
Repeating this construction at every point of
$D$, we see  there is an extension
$$(E,\overline{\nabla})$$
where $E$ is a holomorphic vector bundle on $Y,$
$ E|_X = F,$
\begin{equation} \label{estensione}
\overline{\nabla}: E \to
E(D)\otimes \omega_Y =
E(D)\otimes \omega_Y .
\end {equation} If $E_{j}$
is the fiber over $P_{j}\in D$, the residues
define maps:
$$\Gamma_{j}=Res_j({\nabla}) : E_j\to E_j,$$
and (up to conjugation)
$\exp(-\pi i\Gamma_j)$ gives the
holonomy of $\nabla$ around  $P_j.$
The eigenvalues $\alpha_{j}$ and
 $-\alpha_{j}$ of $\Gamma_{j}$
are real with $ 0\leq
\alpha_{j}< 1.$ We define the weight
   divisor $$\sum_i \alpha _j P_j.$$  To
define the parabolic structure
(see \ref{parabolic}) we set, for
$\alpha_{j}>0,$
 $$W_{i}=
\ker (\Gamma_i - \alpha_i id).$$
Since $F$ is a $SU(2)$ flat vector bundle,
$\det F$
is trivial  and
$\wedge^{2}\nabla=d$ is
the trivial connection.
Arguing as in
\ref{bandet} we see
$\wedge^{2}\overline\nabla$ extends
$\wedge^{2}\nabla$ and gives a
regular connections on
$\det E$, that
is, $\{\det E,\wedge^{2}\overline\nabla\}$
 is a flat unitary line-bundle
with trivial monodromy: $\det E=\cO_{Y}.$

\begin{definition} \label{parabundle}
When $\{E,\mathcal P\}=
\{ \sum \alpha_{i}P_{i}, W_i \}$
is the special parabolic
vector bundle  associated to $\{F,\nabla\}$
we write
\begin{equation}
\label{funtore}
 \tau (\{F,\nabla\}) = \{E, \mathcal P\}.
\end{equation}
\end{definition}

It follows that $\{E,\mathcal P\}
=\tau (\{F,\nabla\})$ is semistable,
and stable if the representation is irreducible.
 This gives an
almost invertible functor:

\begin{theorem} \label{parabolico}
The correspondence
$$\{F,\nabla\} \mapsto \tau (\{F,\nabla\}) = \{E, \mathcal P\}
$$ defines a one-to one
functor between special stable
 parabolic bundles and irreducible
$SU(2)$-bundles. Moreover, if
 $\{E, \mathcal P\}$ is
 semistable parabolic, there
is an $SU(2)$-bundle
$\{F,\nabla\}$ such that $\tau
(\{F,\nabla \})= \{E, \mathcal P\}.$
\end{theorem}
\begin{proof} See \cite{methases}
and the appendix of \cite{belkale}.
\end{proof}
\bigskip

We recall now a basic result we wiil need in section $4$.

\begin{lemma} \label{propria}
 Let $s$ and $t$ be linearly independent
  meromorphic global sections of an
 irreducible $SU(2)$-bundle $\{E,\co\}$
 (i.e.  the associated $\{E,\cP\}$ is
 stable).  If $s\wedge t =0$, then
 $ \co s\wedge \co t\neq 0.$
\end{lemma}
\begin{proof}
The proof is standard.
We assume
$ \co s\wedge \co t =0$ and we show that
$\{E,\co\}$ is reducible.  First, since
$s\wedge t =0,$ we find a non-constant
meromorphic function $g$ such that
$t=gs.$ Then one has :  $$0=\co s
\wedge \co t = \co s \wedge \co (g\cdot s)=
(dg)s\wedge \co s+ g\cdot \co s
\wedge \co s= dg\cdot (s\wedge\co s),$$ and hence
  $$ s\wedge\co s=0$$
  because $g$ is non constant.
   We obtain
  that the sub-bundle generated by $s$
  is holonomy-invariant. Then $\{E,\co\}$ is reducible.
 \end{proof}

\subsection{ Algebraic cmc-$1$ surfaces
in $\mathcal H^3$ }

Fix a special semistable parabolic
vector bundle $\{E, \mathcal
P\}$ where
$ \mathcal P = \{ D(\cP), W_i \}$,
$\ D(\cP) =\sum
\alpha_{i}P_{i}$, and  $D_{\cP}= supp \, (D(\cP)=
\sum_{i}P_{i}$. Let
 $\co:E \to E(D)\otimes \omega_Y$ be
  the singular connection
 (see \ref {parabolico})
whose
restriction to $Y'=Y\setminus D(\cP)$
gives an $SU(2)$ flat
connection $\nabla$ on $F=E|_{Y'}.$
We recall that $\det E= \cO_{Y}.$

\begin{definition} Two meromorphic
sections $(s,t)$ of $E$ define
a special meromorphic frame of $E$
if $s\wedge t\equiv 1.$
\end{definition}

The meromorphic special frame $(s,t)$ of
 $E$ provides a holomorphic section of
$E^{2}\otimes \cO_{Y}(S)$ for a suitable
 divisor $S.$   We can write:
$$1= s\wedge t\in H^{0}(Y,
\det (E^{2} \otimes \cO_{Y}(S)))=
H^{0}(Y, \cO_{Y}(2S)).$$

Let $L\equiv O_{Y}(S) $ be the
line bundle associated to $S$
and $\sigma \in H^0(Y,L) $ be the section
corresponding to $1$. This means
that the zero divisor of $\sigma:
\cO_{Y} \to L$ is $S.$ Set $v=\sigma\cdot s$
and $w=\sigma\cdot t\in
H^0(Y,E(L))$. The condition
 $s\wedge t =1$ becomes:
\begin {equation}
 \label{homosp} v\wedge w =\sigma^2\in
 H^0(Y,L^2);\ \sigma\neq 0.
  \end{equation}
If conversely  $L$ is a fixed
line bundle and  $v$ and $w$
in $H^{0}(Y,E(L))$ such that $v\wedge
w=\sigma^{2}\in H^0(Y,L^{2})$ where
$0\neq \sigma \in H^0(Y,L)$ then
$s=\displaystyle\frac{v}{\sigma}$ and
$t=\displaystyle\frac{w}{\sigma}$
define a special meromorphic frame
of $E.$
\begin{definition}\label{specialL}
 A couple of sections $(v, w)$ in
$H^{0}(Y,E(L))$ defines a special
frame of $E(L)$ if $0\neq v\wedge
w=\sigma^{2}$, $\sigma\in H^{0}(Y,L).$
\end{definition}
Let $D$ be a simple divisor and
 $X=Y\setminus D.$

\begin{definition}
A $h$-liftable map $f: X\to \cH^3$
 will be said of finite type
  if it is associated to a special
 meromorphic frame $(s,t)$ of a
 parabolic bundle $\{E,\cP\}.$
 \end{definition}
To a special
meromorphic frame $(s,t)$, we have associated in
 (\ref{quadratic})
 a meromorphic quadratic
differential
 $ \Omega (s,t) = \overline
\nabla s \wedge \overline \nabla t$
 on $Y$.
  We have
 \begin{lemma}
 Let $(s,t)$
 be a special meromorphic
 frame with poles on $S$. Then $ \Omega (s,t)$
 has poles on $2S+D_{\cP},$
  that is:
 $  \Omega (s,t)\in H^0(Y,
 \omega^{2}_{Y}(2S+D_{\cP})).
 $
 \label{glopolo}
 \end{lemma}
 \begin{proof} It follows
 from lemma \ref{polo}.
 \end{proof}
The meaning of the above lemma is the
 following:
 the condition $\det(s,t)\equiv 1$
 forces the quadratic differential
 $ \Omega (s,t)$
to have slightly milder singularities.
The homogeneous form of the above
 quadratic differential is
\begin{equation}
 \Theta (v,w) = \sigma^{2} \cdot \Omega
(\frac{v}{\sigma},\frac{w}{\sigma}) \in
H^0(Y,\omega^{2}_{Y}\otimes L^{2}(D_{\cP}))
\label{homqua},
\end{equation}
$v$ and $w$ in $H^{0}(Y,E(L))$ and
 $v\wedge w= \sigma^{2}.$
We seek for the case
 $\Omega(s,t)=0,$  i.e.  $\Theta (v,w)=0$
 (see \ref {bryant}).
 \begin {definition} With
 the above notation
  \begin{enumerate}
  \item a meromorphic special frame $(s,t)$ of
   $E$ will be said a
  meromorphic $\mathcal B$-frame of $E$ if
  $\Omega(s,t)=0;$
  \item a special frame $(v,w)$ of $E(L)$
  will be called a $\mathcal B$-frame
   of $E(L)$ if $\Theta (v,w)=0;$
  \item a $h$-liftable non-constant map
  $f: X\to \cH^3$ associated to a
  meromorphic $\mathcal B$-frame will
   be called an algebraic cmc-$1$
  curve of  hyperbolic space.
  \end{enumerate}
\label{alg}
\end{definition}
\begin{lemma}  Let $(s,t)$ be a
meromorphic $\mathcal B$-frame of
$E$ and $S$ be its polar divisor.
Then $D= supp \,(S)$ contains the
support of the parabolic divisor:
 $D-D_{\cP}\geq 0$ i.e.  $D-D_{\cP}$
 is effective.\end
{lemma}
\begin{proof} Use  formula  \ref{stareq}.
\end{proof}


\section {Existence results}
 We fix a compact Riemann surface
 $Y$ of genus $g$ and
a semistable
parabolic bundle $\{E, \cP \}$ on $Y.$
The singular connection
 on $E$  will be
 $\overline\nabla,$  its
support divisor 
$D_{\cP}$  and    $d_{\mathcal P}=\deg D_{\cP}$ (\ref{parabolic}).

\subsection {The variety of special frames}

Let $L$ be a holomorphic line bundle on
$Y$ of degree $d.$

\begin{hyp} \label{hyp}
From now on
we assume  $d\geq \frac{1}{2}d_{\mathcal P} +2g.$
 Riemann-Roch
and \ref{vanishing}
 imply:
\begin{enumerate}
\item[a)]  $ h^0(L) =d-(g-1)$;
  $h^0(L^{2}) =2d-(g-1).$
 \item[b)] $h^1(E(L))=0$, that is, $h^0(E(L))= 2d-2(g-1).$
  \item[c)] $E(L)$ is generated by its global sections.
\end{enumerate}
 \end {hyp}
 \noindent We define the basic determinant
map $\phi: H^0(Y,E(L))\times
H^0(Y,E(L))\fr H^{0}(Y,L^{2})$
\begin{equation}
 \label{det}
  \phi(\omega_{1}, \omega_{2})=
  \omega_{1}\wedge \omega_{2}.
 \end{equation}

Let $$\widetilde Q= \{ \vartheta
 \in H^{0}(Y,L^{2}): \vartheta = \sigma^{2} \
: \sigma\in H^{0}(X,L)),\ \sigma \neq 0\}.$$
Let $\mathbb P= \mathbb P( H^{0}(X,L^{2}))$
 be the projective
 space of $ H^{0}(X,L^{2}).$
 The locus $\widetilde Q$
 is a homogenous cone, we
let $\bQ$ be its associated  projective locus:
 \begin{equation}
 \label{pqua}
 \bQ = \{( \vartheta): \vartheta\in \widetilde Q\}.
 \end{equation}
 The Veronese embedding
 $(\sigma)\to (\sigma^{2})$ identifies
 $\bQ$ and $\bP(H^{0}(X,L)),$
 the projective space associated to
 $H^{0}(X,L)$. In particular
 $\dim \bQ= d-g.$ We recall that
 $(\omega_1, \omega_2)$ gives
 raise to a special frame if and only if $
 \phi(\omega_{1}, \omega_{2})\in \widetilde Q.$
 Let  $\mathbb G=
 \mathbb G(2, H^0(Y,E(L))$ be
 the grassmannian of $2-$planes in
 $H^0(Y,E(L)).$ We have $\dim \mathbb G= 4(d-g).$
 \begin{definition}\label{sp}
 The locus $$\cF = \phi^{-1}(\widetilde Q)$$
 will be called the locus
  of special frames of $E(L).$
 A plane generated by a special frame,
  $\Pi= span (\omega_{1}, \omega_{2}):$
  $\phi(\omega_{1}, \omega_{2})\in \widetilde Q ,$
  will be called special. Therefore the locus
 $$
  \cG =\{ \Pi\in \mathbb G: \Pi\ {\rm is \ special}\}
 $$
  will be called the special planes locus.
 \end{definition}
 We remark that any basis of a special
 plane gives a special frame of $E(L).$
We have then a natural fibration
$\psi: \cF \to \cG$  defined by:
 $$
 \psi(\omega_{1}, \omega_{2})= span(\omega_{1}, \omega_{2}).
 $$
Finally we define the
  mapping: $\varrho: \cG \fr \bQ\subset \bP$ by
 \begin{equation}
  \varrho (span(\omega_{1}, \omega_{2}))= (\phi(\omega_{1},
 \omega_{2}))= (\omega_{1}\wedge \omega_{2}).
 \label{grp}
 \end{equation}

Let
$\Upsilon\subset \bG \times \bQ$ be the closure
 of the graph $\Upsilon_{\varrho}$ of
$\varrho.$ That is  :
\begin{equation}
 \label {chiusura}
\Upsilon= \overline{ \{ (\Pi,\varrho(\Pi)):
\Pi\in \cG \}}=\overline
{\Upsilon}_{\varrho}.
\end{equation}
Set $T=\bG \times \bQ.$
We recall that using the Veronese embedding
we have identified $\bQ$ and
$\mathbb P(H^{0}(Y,L))$.
Let $q_{1}:T\to \bG$ and
$q_{2}:T\to \bQ$ be the projections.
We still
 denote by $q_{1} :\Upsilon \to \bG$ and
$q_{2} :\Upsilon \to \bQ$ the
induced projection maps.
We can identify $\cG$ with
 $\Upsilon_{\varrho},$ and
 $q_{1}|_{\Upsilon_{\varrho}}$ with $\varrho.$
 We have $q(\Upsilon)=\overline{\cG}$
 the closure of $\cG.$
 Set
$\Upsilon_{0}= \Upsilon\setminus \Upsilon_{\varrho},$ and
$\cG_{0}=\overline{\cG}-\cG.$

\begin{proposition}
Let $\Pi\in \cG_{0}$ and $(v,w)$
be any basis of $\Pi.$
Then : $v\wedge w=0.$
\label{zero}
\end{proposition}
\begin{proof} If $\Pi=span(v,w)\in \Upsilon$ then its determinant
 $v\wedge w = \sigma^{2}$
 is a square where  $  \sigma \in H^{0}(X,L);$
  $\sigma \neq 0$  if and only if
 $\Pi \in \cG.$
 \end{proof}


\subsection{Existence of special frames}

We show under the hypothesis
\ref{hyp} that the locus of
special frames of $E(L)$ is not empty. Let $v\in
H^{0}(Y,E(L))$ be a fixed section.
 Let $\phi_{v}: H^{0}(Y,E(L)) \fr
H^{0}(X,L^{2})$ be the linear map :
 $$
 \phi_{v} (\omega)= \phi(v,\omega)= v \wedge \omega.
 $$
  We have a well-known:
  \begin{lemma}\label{glob}
  If $v$ is a general section then $\dim
 \ker(\phi_{v})= 1,$ i.e. $\ker (\phi_{v})$
  is generated by $v.$
 \end{lemma}
 \begin{proof}
 Since $E(L)$ is generated by global
  section we can find
 $v\in H^{0}(Y,E(L))$ without zeros.
  Assume $v \wedge \omega=0,$
then $\omega= f\cdot v$ where $f$ is
 a meromorphic function on $Y.$ Since
 $v$ has no
zeros, then $f$ has no poles,
and so it is constant.
 \end{proof}
 Denote by
$\Phi_{v} =\{(a)\in \cP, \ \ a= v\wedge \omega \}.$
It follows
from lemma \ref{glob} that $\Phi_{v}$
is a projective space of dimension
$h^{0}(E(L))-2=2d-2g$, that is, $\Phi_{v}$
has codimension $g$ in $\bP.$
Since $\bQ$ is a projective subvariety
of codimension $g$ of $\bP,$  we have proved:

\begin{proposition}
\label{phi}
 The locus $ \bQ_{v}= \bQ \cap \Phi_{v}$
 is projective of codimension
 $e\leq d+g.$
 That is, $\dim \bQ_{v} \geq d-2g$
 and, in particular, (by
 \ref{hyp}) it is not empty.
\end{proposition}

A dimension count gives:

\begin{corollary} \label{framedim}
  The loci $\cF$ and $\cG$ are
 not empty and moreover:
 $$\dim \cF \geq 3d -(4g-4) \ , \ \dim\cG\geq 3d - 4g.$$
 \end{corollary}
 \begin{proof}
 The algebraic locus $\cF\subset H^{0}(Y,E(L))^2$
 is not empty by
 \ref{phi}.
  Let $p_{i}: H^{0}(Y,E(L))^{2}\fr H^{0}(Y,E(L)),$
   $i=1,\ 2$
 be the projections.
 Then $ p_{1}(\cF)$ contains the general point $v$
 and so an open Zariski
 set.  Finally
  $p_{2}(p_{1}^{-1}(v) \cup \cF) $
 is the locus
 $$
 \Gamma_{v}=\{ w\in H^{0}(Y,E(L)):
 (v\wedge w)=(\phi_{v}(w))
 \in \bQ \}=\phi_{v}^{-1}(\bQ_{v}) .$$
Since (see \ref{phi})
 the codimension
$\bQ_{v}$ is $d$ we have dimension
$\dim \Gamma_{v}= d-2g+2.$
 We finally obtain
$\dim \cF = \dim \Gamma_{v}+ h^{0}(E(L))=
d-2g+2+2d-2g+2=3d-4g+4$.
 \end{proof}

\subsection {The variety of $\mathcal B$-frames }
Now we will
study the loci defined by the Bryant
 condition (see \ref{bryant} ).  If $(v/\sigma,w/\sigma)$
is a special frame we
defined (see \ref{homqua}) :
$$\Theta(v,w)= \sigma^{2}\cdot
 \co (\frac{v}{\sigma}) \wedge
 \co(\frac{w}{\sigma}) $$
 where $v\wedge w= \sigma^{2}$ and
  $0\neq \sigma \in H^{0}(Y,L)).$
 We also recall (see  \ref{sp}) that the locus of the
  special frames is
  $\cF= \{(v,w)\in H^{0}(E(L))\times H^{0}(E(L)) :
  v\wedge w =\sigma^{2},
   \sigma\in  H^{0}(L), \sigma\neq 0 \}.$

  Since
 $\Theta(v,w)$ is invariant
 under the changing of sign in $\sigma,$
 it defines a map:
 \begin{equation}
  \label{quamap}
 \Theta: \cF \fr H^{0}(Y,\omega_{Y}^{2}(D_{\cP})).
 \end{equation}
 The equation $\Theta(v,w)=0$
 holds (see  \ref{alg})  for the $\mathcal
 B$-frame $(v,w)$ of $E(L).$ We also defined
 $\psi: \cF \to \bG,$
$\psi(\omega_{1},\omega_{2})=
span(\omega_{1},\omega_{2}). $
Accordingly we give the following:
\begin{definition}
A special plane
$\Pi \in \cG \subset \bG$
will be called a $\cB$
plane if it is generated by a
$\cB$-frame.
We set:
\begin{enumerate}
\item[a)] $\cB= \Theta^{-1}(0)\subset \cF \equiv$
the locus of $\mathcal B$-frames.
\item [b)] $ \cM=\psi(\cB)\subset \bG \equiv$
 the locus of $\cB$-planes.
\end{enumerate}
\end{definition}
\noindent We remark that any basis of a
$\cB$-plane is a $\cB$-frame of $E(L).$
\medskip
\begin{remark}
The locus of the special planes
 $\cG$ is not compact
 in general, and moreover (see
\ref{grp})
$\varrho: \cG \to \bQ $ does not extend
to the closure ${\overline
\cG} \subset \bG.$
 We notice however that if
 the condition $v\wedge w=
\sigma^2$ is dropped
formula (\ref{homqua})
is  well defined for $\sigma\neq0$.
\end{remark}
We defined $T=\bG \times \bQ$
with projections
$q_{1}:T\to \bG$ and
$q_{2}:T\to \bQ$.
Let $S_{\bG}$ be the
tautological rank $2$ vector
bundle of $\bG$, that is:
 $$
 S_{\bG}=  \{(\Pi,v)\in \bG\times H^0(Y,E(L)): v\in \Pi \}.
 $$
The line bundle
$\bigwedge^{2} S_{\bG}$
 generates the Picard group of
$\bG$ (\cite {griffithsha} chapter 2)
and
its dual
is ample.
Let $\cO_{\bQ}(-1)$ be
the tautological line bundle of
$\bQ=\bP(H^{0}(X,L))$. The
pull-backs:
$M=q_{1}^{\ast}(\bigwedge^{2} S_{\bG})$ and
$N=q_{2}^{\ast} (\cO_{\bQ}(-1))$
 are line bundles on $T.$ Define
 the trivial bundles
 $$
  \bV=
  H^0(Y,\omega^{2}_{Y}L^4(\cD_{\cP}))\times T,\  \bW=
  H^0(Y,\omega^{2}_{Y}L^2(\cD_{\cP}))\times T.$$

  Given $v$ and $w$ in $H^{0}(Y,E(L))$
 the formula
 \begin{equation}
 \Psi(v,w,\sigma)= \sigma^{4}
\co \frac{v}{\sigma} \wedge \co
\frac{w}{\sigma}
 \label{Psi0}
 \end{equation}
 is well defined if
 $0\neq \sigma \in H^{0}(X,L)$
  and $\Psi(v,w,\sigma)\in \bV .$
 Letting $\lambda
   \in \mathbb C, $ $A\in GL(2,\mathbb C),$
   we observe the homogeneity:
 $$ \Psi(Av,Aw,\lambda\sigma) =
  \lambda^{4} \sigma^{4}
\co \frac{Av}{ \lambda\sigma} \wedge \co
\frac{Aw}{ \lambda\sigma}=
 \det A \lambda^{2} \Psi(v, w,\sigma).
  $$
This defines a bundle map:
 \begin{equation}
\Psi: M\otimes N^{2}\to \bV.
 \label{Psi00}
 \end{equation}
Now we study the restriction of
$\Psi$ to $\Upsilon.$
 First we see
by  \ref{glopolo} that on
$  \Upsilon_{\varrho}$ we have
$$ \Psi(\sigma, v,w)= \sigma^{2} \Theta(v, w)$$  where
$\Theta(v,w) \in  H^0(Y,\omega^{2}_{Y}L^2(\cD_{\cP})).$
Then the locus   $$\Gamma= \{(span(v,w),(\sigma)) \in T :
\Psi(\sigma,v,w)=
\sigma^{2}\eta, \ \eta \in
H^0(Y,\omega^{2}_{Y}L^2(\cD_{\cP}))\}$$
is closed and  $\Upsilon_{\varrho}\subset \Gamma$,
so that
 $\Upsilon \subset \Gamma.$
In fact if $\Pi=span(v,w)\in \Upsilon$ we have:
$$
\Psi(v,w,\sigma)\in \sigma^{2}\bW\subset \bV.
$$
This observation is essentially the content of
 \ref{glopolo}.
  We will denote with $\bW_{\Gamma}$ and $M_{\Gamma}$
  the restriction of $\bW$  and $M$
 to $\Gamma .$
Let, for $\sigma\neq 0,$
$$\overline \Psi (v,w,\sigma)=
\sigma^{-2}\ \Psi(v,w,\sigma).$$
We have the homogeneity:
$$
\overline \Psi(Av, Aw,\lambda \sigma) =
\det(A) \overline \Psi(v,w,\sigma)
$$

\noindent Therefore the restriction of $\Psi$ to
$ \Gamma$ defines a bundle map
$\overline\Lambda : M_{\Gamma}\to \bW_{\Gamma}$, and
this bundle map is well defined on
$\Upsilon\subset \Gamma.$ Let
$M|_{\Upsilon}=\widetilde{M}$
 be the restriction of $M$ to
$\Upsilon$ and denote for sake of notation  by
$\bW_{\Upsilon}$ the restriction of  $\bW$ to $\Upsilon.$
It follows that the restriction of  $\Psi$ to
$\Upsilon$ defines a
vector bundle map:
\begin{equation}
\label{Psi}
\Lambda: \widetilde{M}\to \bW_{\Upsilon}.
\end{equation}

In the case of irreducible monodromy
the boundary loci $\cG_{0}$
and $\Upsilon_{0}$
do not intersect the zeroes of $\Lambda.$
In fact we have the following:
\begin{lemma}
Assume that $E$ is (parabolic)
stable. Then $ Z_{\Lambda}\cap
\Upsilon_{0}=\emptyset.$
In particular the projection
$q_{1}=\varrho:\Upsilon_{\rho} \to \bG$
 defines an isomorphism
 between
$Z_{\Lambda}$ and the locus
$\cM$ of $\cB$ planes.
\label{interno}
\end{lemma}
\begin{proof}
It is the content of lemma \ref{propria}.
\end{proof}
Now we state our existence result:
\begin{theorem}\label{esistenza}
Assume that $\{E,\co \},$
is associated to an irreducible $SU(2)$-flat
representation, i.e.,
$\{E, \cP\}$ is stable.
We set $\deg D_{\cP}= d_{\cP}.$
Let $L$ be a holomorphic
line-bundle of degree $\deg L=d.$
Suppose $$ d\geq
7g-3+d_{\cP}.$$Then the locus
$\cM$ of $\cB$-planes
of $E(L)$ is a nonempty algebraic
set of dimension $a\geq
d-7(g-1)-d_{\cP}-4.$
\end{theorem}
\begin{proof}
By \ref{interno} we have to show
that $Z_{\Lambda}$
is not empty
of dimension $a\geq d-7(g-1)-d_{\cP}.$
We  consider
 $\Lambda\in H^{0}(\bW_{\Upsilon} \otimes (\widetilde M)^{-1}).$
 Since by \ref{framedim}
$\dim \Upsilon=\dim\cG\geq 3d-4g$
it is enough
 to show the top Chern class $c_{top}$
 of  $\bW_{\Upsilon} \otimes (\widetilde M)^{-1}$
 is not zero. Since $\bW_{\Upsilon}$
 is trivial we have
 $$
   c_{top}=  c_{1}({\widetilde
 M}^{-1})^{r},
 $$
 where $r = 2d+3g-3+d_{\cP}$
 is the rank of $\bW_{\Upsilon}$.
 Let $i: \Upsilon\hookrightarrow\bG$
 be the inclusion.
It follows then
 $$
 c_{top}=i^{\ast}(q_{1}^{\ast}(h)^{r}),
  $$
where $ h= c_{1}(\wedge^{2}S_{\bG})^{\ast}$
 is positive since
 $\wedge^{2}S_{\bG}^{\ast}$ is ample.
 We compute the cohomology class
 $\mu= q_{1 \ast} i_{\ast} ( c_{top})$ of $\bG.$
First we obtain:
 $$
 i_{\ast} c_{top}=
 i_{\ast}i^{\ast}(q_{1}^{\ast}(h)^{r})=
 q_{1}^{\ast}(h^{r})\cdot [\Upsilon],
  $$
 where $ [\Upsilon]$
 denotes the class of the variety
 $ \Upsilon\subset T.$
 Similarly we let $[{\overline \cG}]$
 be the class of
 ${\overline \cG}$ in $G.$
 Then we get
 $$\mu= q_{1 \ast}
 q_{1}^{\ast}(h^{r})\cdot [\Upsilon]=
 h^{r} \cdot q_{1 \ast}[\Upsilon]=
 h^{r}\cdot [{\overline \cG}] \neq 0$$
since $h$ is positive and the dimension of
${\overline \cG}$ is larger than $r$.
It follows
 that  $c_{top} \neq 0,$
which proves the theorem.
\end{proof}

\noindent {\em Proof of the Existence Theorem}.\label{ex}
One observes that \ref{esistenza} provides
$\mathcal B$-meromorphic frames of $E$ in the case of
any irreducible monodromy on $Y-D.$
If $(s,t)$ is such a $\cB$-frame, let $P_{s}$
and $P_{t}$  be the polar divisors  of $s$ and $t$ respectively.
Let $Z=\{z\in Y :\nabla s (z)=\nabla t (z)=0\}$
the branch points divisor.
Let $H=P_{s}\cup P_{t}\cup Z$ be the union.
We  define a cmc-$1$ immersion (see \ref{global})
$f:Y\setminus ( D\cup H )\to \mathcal H^{3}.$
By construction the monodromy of $f$
is the monodromy of $E$.


\section {Higgs fields and the period problem}

In this section we use the special
frames to obtain a standard form for the flat
connection.
 Then the period problem is translated into a holonomy
problem.
The language of Higgs fields (see
 \cite {hitchin} and
 \cite{simpson})
allows,
in particular, to discuss the
conditiont given in \cite{ktuy}.


\subsection{Higgs fields associated to h-liftable map}

Let $G=SL(n,\bC)$
and let $H$ be a closed subgroup of $G$.
Let  $H_{e}$
 be the connected component of $H$ containing $e$:
 we assume $H_{e}=SU(n)$.
 Let $f: X \to SL(n)/H$ be h-liftable,
 $\{F, \nabla\}$ be the
  associated flat bundle.  The
  associated special frame, $e_{1},\ldots,e_{n}$,
  gives a holomorphic trivialization of $F:$
  $$F\ \equiv _{f} \ \cO_{X}^{n}.$$
  Since $f$ is $h$-liftable, we have
  $\nabla e_{i}\in H^{0}(X, \omega_{X})^{n}$. We write
  $$ \nabla e_{i}= (\omega_{1,i},\ldots \omega_{n,i}).$$
Transposing we obtain a  matrix of $1.0$-forms.
\begin{equation} \label{mhiggs}
 \Theta=^{t}\!\!(\nabla e_{1},\ldots, \nabla
e_{n}) = (\omega_{i,j})
 \end{equation}
 We will identify $\Theta$
 with the associated map
 $\Theta: F\to F\otimes
 \omega_{X},$ that is $\Theta
 \in H^{0}(X, Hom(F,F\otimes \omega_{X})).$
We may write $$\nabla = d+ \Theta.$$
The decomposition $\nabla=D'+D''$ gives
$D'=\partial + \Theta$ and
$D''=\overline{\partial}.$
The condition that $\nabla$ is special translates into
the vanishing of
the trace of $\Theta.$ In fact one has
$$
0= (\wedge^{n}\nabla)(e_{1}\wedge \ldots \wedge
e_{n})=(\sum_{i=1}^{n}\omega_{i,i})e_{1}\wedge \ldots \wedge
e_{n}.
$$
and hence
 $$\sum_{i=1}^{n}\omega_{i,i}=0.$$
It follows that $\Theta$ is a
$sl(n,\bC)$ matrix.
\begin{definition}
The matrix $ \Theta$ in (\ref{mhiggs}) will be called the
Higgs field associated to $f$.
\label{higgs}
\end{definition}
 Let
 $M_{\Theta}\subset SL(n,\bC) $
 be the holonomy group of $\nabla.$
The topological
 closure $K_{\Theta}=\overline{M_{\Theta}}$
   of $M_{\Theta}$
in $SL(n,\bC)$ (by Cartan theorem)
defines a Lie-subgroup. By construction one
has that the connected component of
$K_{\Theta}$ through $e$ is contained in $H=SU(n).$

Conversely given on $F=\cO_{X}^{n}$  a connection
 $ \nabla = d + \Theta$ where $\Theta=(\omega_{i,j})$
is a tensor of type $(1,0)$,
 we have:
  \begin {enumerate}
  \item $ \nabla^{2}=0$ if and only if
    $\overline{\partial}(\Theta)=0$ (i.e.  the
    $\omega_{i,j} $ are holomorphic).
    \item $\nabla$ is  a $SL(n,\bC)$ connection  if and only if
    the trace of $\Theta$ vanishes, i.e.
     $\sum_{i=1}^{n}\omega_{i,i}=0.$
      \end{enumerate}

Given a holomorphic trace-free
Higgs field $\Theta$,
let $\nabla_{\Theta}=d + \Theta$ be
the associated connection. Let $M_{\Theta}$ be
its holonomy group and $K_{\Theta}$ be its closure in $G$.
We give the following:
\begin{definition} \label{pprob}
 ({\bf Non abelian period problem}): We say
that $\Theta$ satisfies the period
problem if the connected component
of $K_{\Theta}$ through $e$ is contained in $SU(n).$
\end{definition}
 We may state:

\begin{proposition} \label{scianti}
Assume $\Theta$ satisfies the period problem of
\ref{pprob}.
Then the standard basis of $\cO_{X}^{n}$ defines a
$h$-liftable map $f: X \to G/H$ where
$H\supset K_{\Theta}$ and $H_{e}=SU(n).$
\label{problemaperiodi}
\end{proposition}
\begin{proof}
For any simply connected open set $U$ of $X$
fix flat unitary basis $ f_{1},\ldots f_{n}$ of $\bC^{n}$
 consider the matrix $A=(a_{i,j})$ such that
 $e_{i}=\sum_{j}a_{i,j}f_{j}$. Then define
 $f(x)=A(x) \cdot A^{\ast}(x) \mod H$ for $x\in U.$
\end{proof}

\begin{remark} It appears difficult
to provide conditions in order to solve the above period problem.
\end{remark}

\subsection {Higgs fields and Maurer-Cartan form}

 We consider again a $ h$-liftable
 map $f:X\to G/H$ where $G=SL(n,\bC)$
 and $H_{e}=SU(n)$ as in the previous
 section. We let $\Theta$ be the Higgs
 field associated to $f$ (see \ref {higgs})
  and $\nabla= d+\Theta$ be
 the associated flat connection on $F= \bC^{n}.$
Let $\{U,z\}$ be a simply
connected coordinate open set.
 Using the canonical frame of $\bC^{n},$
  we set $E=(e_{1},\ldots,e_n).$
  As in proposition  \ref {scianti},
  we have  matrices $A(z)$
  such that
   $$f(x)=A(z(x))\cdot A^{\ast}(z(x))  \mod H$$
   for $x\in U. $
  Moreover $A^{-1}E = R =(f_{1},\ldots f_{n})$
  is a  unitary frame of $F|_{U}$.
  Taking derivative in the equation
  $A\cdot R=E$ we find
  $$ dA \cdot R= \nabla (A\cdot R) =
  \nabla E= \theta \cdot  E=\theta A\cdot R $$
 that is:
 \begin{equation} \label {m-c}
 \Theta = dA \cdot A^{-1} .
 \end{equation}
 Now comparing with \cite [section 3] {ktuy}
 we see that $\Theta$
 is pull-back to $U$  of the
 {\em right invariant}
 Maurer- Cartan form of $G=SL(n,\bC)$,  that is:

\begin{proposition} The Higgs field
$\Theta$ is the pull-back of the right-
invariant Maurer Cartan form to $X.$
\end{proposition}

If we now consider the matrix
of quadratic forms:
$$ \Theta ^{2}= (\beta_{i,j})=
\ \sum_{k=1}^{n} \omega_{i,k}\omega_{k.j} $$
we have the following generalization of the Bryant condition:
\medskip

\noindent
{\bf Kokubu-Takahashi-Umehara-Yamada condition}
(see \cite  [3.3] {ktuy})
\begin{equation}
{\rm trace} (\Theta ^{2}) = \sum_{k=1}^{n} \beta_{k,k}=
\sum_{i,k} \omega_{i,k}\omega_{k,i}  = 0.
\label{ume}
\end{equation}

We specialize formula (\ref {ume}) to
the case $n=2$ to get the Bryant
condition  (\ref {bryant}), that is
$\nabla e_{1}\wedge \nabla e_{2}= 0.$
If we write the $sl(2)$ matrix:
 $$ \Theta= \begin{pmatrix} \alpha & \gamma \\
  \beta & -\alpha
\end{pmatrix}
$$
\noindent we obtain $ \nabla e_{1}=
 \alpha e_{1}+ \beta e_{2}$
and $ \nabla e_{1}=
\gamma e_{1}- \alpha e_{2},$ and hence
$$\nabla e_{1}\wedge \nabla e_{2}=
 -\alpha^{2}- \beta\gamma=\det \Theta.$$
Then the frame satisfies the Bryant
condition if and only if $\det \Theta=0.$
That is
$\Theta$ is nihilpotent:
$\Theta^{2}= \det\Theta {\rm I} $ or equivalently
the  trace of $\Theta^{2}$ vanishes.

\begin{remark}
Cousin's representation of
a (local) minimal surface
in Euclidean space is obtained,  when $ \det\Theta =0,$  taking
$\omega_{1}=\alpha,$
$\omega_{2}= i(\beta+\gamma)$  and \
 $\omega_{3}=i(\beta-\gamma).$

\end{remark}



                                %
                                %
                                %
                                %
                                %
                                %

\noindent {\sc Gian Pietro Pirola}\\
Dipartimento di Matematica, Universit\`a di Pavia\\
via Ferrata 1, 27100 Pavia,  Italia\\
 fax +39 0382 985602\\
{\tt em: gianpietro.pirola@unipv.it}

\end{document}